\documentclass[12pt,final]{amsart}

\usepackage{graphicx}
\usepackage{psfrag}
\usepackage{subfigure}
\usepackage{color}
\usepackage{bm}
\usepackage{epstopdf}

\usepackage{showkeys}

\usepackage{paralist}

\numberwithin{equation}{section} \numberwithin{figure}{section}
\def\n{\noindent}
\def\pt{\partial}

\def\Re{\mathbb{R}}

\def\f#1#2{\frac {#1}{#2}}
\def\f32{\frac 32}

\def\d{\displaystyle}

\def\bga{\begin{array}}
\def\eda{\end{array}}

\def\De{\Delta}

\def\gm{\gamma}

\def\nb{\nabla}

\def\al{\alpha}

\def\iy{\infty}
\def\Om{\Omega}
\def\om{\omega}

\def\sq{\sqrt}

\def\d{\displaystyle}
\def\dfr#1#2{\displaystyle{\frac{#1}{#2}}}

 \newtheorem{thm}{Theorem}[section]

 \theoremstyle{definition}
 
 \theoremstyle{remark}
 \newtheorem{rem}[thm]{Remark}
 
\def\bbf{\mathbf{f}}

\def\bn{\mathbf{n}}

\def\bu{\mathbf{u}}

\def\bx{\mathbf{x}}

\def\bA{\mathbf{A}}

\title{A Hermite WENO reconstruction for fourth order temporal accurate schemes  based on the GRP solver for hyperbolic conservation laws}

\author{Zhifang Du and Jiequan Li}

\address{Zhifang Du: School of Mathematical Sciences, Beijing Normal University, 100875, P. R. China; Email: du@mail.bnu.edu.cn}

\address{Jiequan Li: Laboratory of Computational Physics,  Institute of Applied Physics and Computational Mathematics, Beijing, P. R. China; Email: li\_jiequan@iapcm.ac.cn}\thanks{This research is supported by NSFC with No 11371063.}

\begin{document}
\maketitle
\markboth{}{A HWENO scheme  based on the GRP solver}

\begin{abstract} 
This paper develops a new fifth order accurate  Hermite WENO (HWENO) reconstruction method for hyperbolic conservation schemes in the framework of the two-stage fourth order accurate temporal discretization in  [{\em J. Li and Z. Du, A two-stage fourth order time-accurate discretization {L}ax--{W}endroff type flow solvers, {I}. {H}yperbolic conservation laws, SIAM, J. Sci. Comput., 38 (2016), pp.~A3046--A3069}].  Instead of  computing the first moment of the solution additionally in the conventional HWENO or DG approach,  we can directly take the {\em interface values}, which are already available  in the numerical flux construction using the generalized Riemann problem (GRP) solver, to approximate the first moment. The resulting scheme is fourth order temporal accurate by only invoking the HWENO reconstruction twice so that it becomes more compact.  Numerical experiments show that such compactness   makes significant  impact on the resolution of nonlinear waves. 
\end{abstract} 


\section{Introduction}
In the development of high order accurate schemes for hyperbolic conservation laws,  two families of approaches  play important roles: one belongs to the  method of line  that achieves the temporal accuracy using the Runge-Kutta strategy \cite{Harten-ENO, Shu-Osher-88, Liu, Jiang, Shu-16}; the other is the Lax-Wendroff type approach that adopts the Cauchy-Kowaleveski expansions to design temporal-spatial coupled schemes \cite{L-W, Ben-Artzi-84, Li-1, menshov-1990, Qiu-Dum-05,Toro-98, Xu-93}.   Either family of approaches have their own advantages and disadvantages.  The former has the simplicity in their practical implementation thanks to  exact or approximate Riemann solvers, but the multi-stage 
temporal iteration inevitably causes the enlargement of the size of stencils; the latter can avoid the multi-stage  temporal iteration but have to  repeatedly make the differentiation of governing equations in order to construct high order accurate numerical fluxes.  A recent two-stage fourth order accurate  temporal discretization based on the Lax-Wendroff type solvers \cite{Du-Li-1, xu-li} makes a compromise between these two families of methods:  It just takes  a two-stage iteration for the fourth order accuracy by using second order accurate temporal-spatial coupled Lax-Wendroff flow solvers so half of reconstruction steps can be saved in comparison with the same accurate method and complicated successive  differentiations of governing equations can be avoided, which could be further extended using the multi-derivative Runge-Kutta methods \cite{xu-li, Seal-16, 
Okten-17}.  Moreover, we notice that  the solution values on cell interfaces already available  in the procedure of numerical flux construction, called {\em interface values} in the present paper,  can be used for the reconstruction procedure, thanks to the Lax-Wendroff flow solvers, which motivates us for such a study.\\

We develop a new fifth order accurate Hermite WENO (HWENO) reconstruction in the framework of two-stage fourth order accurate temporal discretization \cite{Du-Li-1}.  The HWENO interpolation adopts two values: the average value of the solution and the corresponding averaged gradient value (the first moment), as usual.  The novelty is that the gradient values are directly approximated using  the interface values when the Lax-Wendroff type flow solvers are used \cite{Ben-Artzi-84, Li-1, Xu-93,Toro-ADER}, which is different from the standard HWENO method in \cite{Qiu-11, Qiu-12, Luo-07}.  Technically,  we can further adjust nonlinear weights during the HWENO reconstruction, just like the WENO-Z method \cite{wenoz}  modifying the classical WENO-JS \cite{Jiang}. In doing so, the resulting scheme is much more compact and   has several distinct features.

\begin{enumerate}

\item[(i)] The scheme just uses half of the reconstruction steps, compared with the standard RK-WENO methods. 

\vspace{0.2cm}

\item[(ii)]  The interface values are already available in the computation of numerical fluxes and no extra efforts  are made on the gradient approximation. 
\vspace{0.2cm}

\item[(iii)]  The interface values are approximated by using the GRP solver and thus they are strong solution values without taking account of possible discontinuities in trouble cells. 
\vspace{0.2cm}

\item[(iv)]  A single HWENO reconstruction  is more compact than the standard WENO reconstruction \cite{Jiang},  as shown in other HWENO schemes \cite{Qiu-11, Qiu-12, Luo-07}. 

\vspace{0.2cm} 

\end{enumerate}

This paper is organized as follows. In Section 2, we quickly review the two-stage method based on the Lax-Wendroff flow solvers and the HWENO reconstruction methods. In Section 3, we show the gradient approximation over each computational cell by using interface values of solutions.  In Section 4, several numerical examples are displayed for the performance of such a HWENO reconstruction, by comparing with the WENO reconstruction with the same numerical flux.  A discussion is made in Section 5.

\vspace{0.2cm}

\section{The two-stage fourth oder method and the Hermite WENO reconstruction}\label{sec:HWENO}

This section serves to present a quick review of the two-stage fourth order method based on the Lax-Wendroff type flow solvers in \cite{Du-Li-1} and the HWENO  reconstruction procedure, originally in \cite{Qiu-11}.  Instead of  independently computing the first moment (the gradient of solution) in \cite{Qiu-11}, we will construct it together with the solution average using the generalized Riemann problem (GRP) solver, which will be described in Setion \ref{sec-moment}.

\subsection{Review of the two-stage fourth-order scheme}\label{sec:time-dis}
%Certainly, we can use the acoustic version of GRP or ADER \cite{Toro-98, Toro-ADER} as no strong waves are involved. Otherwise, as waves involved are very strong, the nonlinear version of GRP or Men'shov's second order solver \cite{menshov-1990} has to be used, as shown in  \cite{thermo-GRP} for theoretical justification and numerical evidence.
The two-stage fourth-order finite volume schemes based on the GRP solver was developed in \cite{Du-Li-1}. Certainly, we can also use Men'shov's modified GRP solver \cite{menshov-1990, menshov-1991} and the ADER solver \cite{Toro-98, Toro-ADER}. Both the acoustic and nonlinear versions of the GRP solver are provided in \cite{Li-1}.

In this subsection, we quickly review this method  by taking  one-dimensional hyperbolic conservation laws,  
\begin{equation}\label{eq:govern-eq}
\begin{array}{l}
  \dfr{\pt \bu}{\pt t} + \dfr{\pt \bbf(\bu)}{\pt x}=0, \ \ \  x\in \mathbb{R},  t > 0,\\[3mm]
  \bu(x,0) =\bu_0(x), \ \ \ \ x\in \Re,
  \end{array}
\end{equation}
where $\bu$ is a vector of conservative variables and $\bbf(\bu)$ is the associated flux function vector.  Given the computational mesh $I_j=(x_{j-\frac 12}, x_{j+\frac 12})$ with the size $ h = x_{j+\frac 12}-x_{j-\frac 12} $ for every $ j $, we write \eqref{eq:govern-eq} in form of the balance law,
\begin{equation}
\dfr{d\bar \bu_j(t)}{dt} =\mathcal{L}_j(\bu):= -\dfr{1}{h} [\bbf(\bu(x_{j+\frac 12},t)) - \bbf(\bu(x_{j-\frac 12},t))], \ \ \ \bar\bu_j(t) =\dfr{1}{h} \int_{I_j} \bu(x,t)dx, 
\label{scheme}
\end{equation} 
where   $\bu(x_{j+\frac 12},t)$ is described in terms of GRP solver \cite{Li-1}.  Then the two-stage approach for \eqref{eq:govern-eq} is summarized as follows.\\

\begin{enumerate}
\item[\bf Step 1.]  With the cell averages $ \bar\bu^n_j $ and interface values $ \hat{\bu}^{n}_{j+\frac 12} $, reconstruct the data at $ t^n $ as a piece-wise polynomial function $\bu(x,t^n) =\bu^n(x)$ by the HWENO interpolation that will be described below, and compute the  corresponding GRP value $(\bu_{j+\frac 12}^n, (\pt\bu/\pt t)_{j+\frac 12}^n)$. 
\vspace{0.2cm} 

\item[\bf Step 2.] Compute the intermediate cell averages $\bar{\bu}^{n+\frac 12}(x)$ and the interface values $\hat{\bu}^{n+\frac 12}_{j+\frac 12}$ at $t^{n+\frac 12}=t^n+\frac k2 $ using the formulae,
\begin{equation}
\begin{array}{l}
\bar\bu_j^{n+\frac 12} =\bar\bu_j^n -\dfr{ k }{2h}[\bbf_{j+\frac 12}^{*}-\bbf_{j-\frac 12}^{*}],  \\[3mm] 
\d \bbf_{j+\frac 12}^{*} = \bbf(\bu_{j+\frac 12}^n) + \dfr{ k }{4} \dfr{\pt\bbf}{\pt\bu}({\bu}^{n}_{j+\frac 12})\Big(\dfr{\pt\bu}{\pt t}\Big)^{n}_{j+\frac 12},\\[3mm]
\d \hat{\bu}^{n+\frac 12}_{j+\frac 12} = \bu_{j+\frac 12}^n + \frac {k}{2}\Big(\dfr{\pt\bu}{\pt t}\Big)_{j+\frac 12}^n.
%\bbf(\bu(x_{j+\frac 12},t^n+\frac 14 k )), 
% \d  \bu(x_{j+\frac 12},t^n+\frac 14 k ):= \bu_{j+\frac 12}^n +\frac{ k }{4}\Big(\frac{\pt\bu}{\pt t}\Big)_{j+\frac 12}^n.
\end{array}
\label{eq:1d-inter}
\end{equation} 
where $  k  $ is the time step size and $\dfr{\pt\bbf}{\pt\bu}$ is the Jacobian of $\bbf(\bu)$.  At this intermediate stage, the HWENO interpolation is carried out  once again with the values $ (  \bar\bu^{n+\frac 12}_j,  \hat{\bu}^{n+\frac 12}_{j+\frac 12}) $,   to construct a piecewise polynomial $\bu^{n+\frac 12}(x)$ and  find the GRP value $(\bu_{j+\frac 12}^{n+\frac 12}, (\pt\bu/\pt t)_{j+\frac 12}^{n+\frac 12})$, as done in Step 1. 
\vspace{0.2cm} 

\item[\bf Step 3.]  Advance the solution to the next time level $t^{n+1}=t^n+k$ by
\begin{equation}
\begin{array}{l}
\bar \bu^{n+1}_j =\bar \bu_j^n -\dfr{k}{h} [\bbf_{j+\frac 12}^{4th} -\bbf_{j-\frac 12}^{4th}],\\[3mm]
\d \bbf_{j+\frac 12}^{4th}=\bbf(\bu_{j+\frac 12}^n) +\dfr{k}{2} \Big[\dfr 13 \Big(\dfr{\pt \bbf}{\pt t}\Big)^{n}_{j+\frac 12} +\dfr 23 \Big(\dfr{\pt \bbf}{\pt t}\Big)^{n+\frac 12}_{j+\frac 12}\Big], \\[3mm] 
%\d \bbf_{j+\frac 12}^{4th}=\bbf(\bu_{j+\frac 12}^n) +\dfr{k}{2} \Big[\Big.\dfr 13 \dfr{\pt  \bbf(\bu)}{\pt t} \Big|_{(x_{j+\frac 12}, t^n) } +\dfr 23 \Big. \dfr{\pt  \bbf(\bu)}{\pt t}\Big|_{ (x_{j+\frac 12}, t^*)}\Big], \\[3mm] 
\d \hat{\bu}_{j+\frac 12}^{n+1} = \bu_{j+\frac 12}^n + k\Big(\dfr{\pt\bu}{\pt t}\Big)_{j+\frac 12}^{n+\frac 12},
\end{array}
\label{eq:1d-final}
\end{equation} 
where the notations are 
\begin{equation*}
  \Big(\dfr{\pt \bbf}{\pt t}\Big)^{n}_{j+\frac 12} = \dfr{\pt \bbf}{\pt \bu}(\bu^{n}_{j+\frac 12})\Big(\dfr{\pt \bu}{\pt t}\Big)^{n}_{j+\frac 12}, \ \ \ 
  \Big(\dfr{\pt \bbf}{\pt t}\Big)^{n+\frac 12}_{j+\frac 12} = \dfr{\pt \bbf}{\pt \bu}(\bu^{n+\frac 12}_{j+\frac 12})\Big(\dfr{\pt \bu}{\pt t}\Big)^{n+\frac 12}_{j+\frac 12}.
\end{equation*}
\end{enumerate}

The same procedure can be applied for two-dimensional cases. The details can be referred to \cite{Du-Li-1}. 

 \begin{rem}
 
 For this two-stage fourth order temporal accuracy  method, the HWENO reconstruction is invoked only twice instead of four times  from $t^n$ to $t^{n+1}=t^n+k$. In addition to save the computational cost, the size of stencils is automatically decreased.

\end{rem} 
\vspace{2mm}

\subsection{A Hermite WENO Interpolation}\label{sec:HWENO}
In the approximation to a given function, there are two types of  typical polynomial interpolations: the Lagrangian interpolation and the Hermite interpolation \cite{Num-book}. The former only uses the grid point values of the function, while the latter uses both the function value and its derivative value at the same grid point. It turns out that  the latter uses half of  grid points  to derive the the polynomial approximations of the same degree to the given function compared to the former interpolation.  So it is significant to develop  a HWENO reconstruction technology for high order accurate schemes of hyperbolic conservation laws, as done in  \cite{Qiu-11}, so that the resulting scheme becomes more compact even though the same flux function approximation is adopted.
\\
 
 Let's summarize the HWENO reconstruction in \cite{Qiu-11} for hyperbolic conservation laws in the finite volume framework although almost the same approach  can be applied in the discontinuous Galerkin (DG) framework too \cite{Luo-07}.   Given the average $\bar \bu_j$ and the derivative $\De \bu_j$ of the function  over the cell $I_j$, 
\begin{equation}\label{eq:0-th}
  \bar{\bu}_j = \dfr 1h \int_{I_j}\bu(x,t)dx,\ \ \ \   \De \bu_j = \dfr 1h \int_{I_j}\dfr{\pt \bu}{\pt x}(x,t)dx,
\end{equation}
we want to construct a polynomial  $p(x)$  such that $\bu_{j+\frac 12,-}:=p(x_{j+\frac 12})$ approximates the left limiting value of $ \bu(\cdot,t) $ at $x=x_{j+\frac 12}$.\\

We choose three stencils
\begin{equation}\label{eq:stencils}
  S^{(-1)} = I_{j-1}\cup I_j, \ \ S^{(0)} = I_{j-1}\cup I_j \cup I_{j+1}, \ \ S^{(1)} = I_j \cup I_{j+1}.
\end{equation}
On stencil $ S^{(0)} $, $ \bar{\bu}_{j-1} $, $ \bar{\bu}_{j} $ and $ \bar\bu_{j+1} $ are used to construct a polynomial $ p^{(0)} $  for the interpolation. Hence at $ x_{j+\frac 12} $, we have
\begin{equation}\label{eq:s-0}
  \bu^{(0)}_{j+\frac 12, -} := p^{(0)}(x_{j+\frac 12}) = -\dfr 16 \bar{\bu}_{j-1} + \dfr{5}{6} \bar{\bu}_{j} + \dfr {1}{3} \bar\bu_{j+1}.
\end{equation}
Similarly, $ p^{(-1)} $ and $ p^{(1)} $ are constructed by using $ \bar{\bu}_j $, $ \bar{\bu}_{j-1} $, $ \De \bu_{j-1} $ on $ S^{(-1)} $ and  by using $ \bar{\bu}_{j} $, $ \bar{\bu}_{j+1} $, $ \De\bu_{j+1} $ on $ S^{(1)} $, respectively, 
\begin{equation}\label{eq:s-n1-p1}
\begin{array}{l}
  \bu^{(-1)}_{j+\frac 12, -} := p^{(-1)}(x_{j+\frac 12}) = -\dfr 76 \bar{\bu}_{j-1} + \dfr{13}{6} \bar{\bu}_{j} - \dfr {2h}{3} \De \bu_{j-1},\\[3mm]
  \bu^{(1)}_{j+\frac 12, -} := p^{(1)}(x_{j+\frac 12}) = \dfr 16 \bar{\bu}_{j} + \dfr{5}{6} \bar{\bu}_{j+1} - \dfr {h}{3} \De\bu_{j+1}.
\end{array}
\end{equation}
If the solution is smooth on the large stencil $ I_{-1}\cup I_0 \cup I_1 $, we have
\begin{equation}\label{eq:s-total}
  \tilde\bu_{j+\frac 12,-} = \dfr{1}{120}(-23\bar{\bu}_{j-1}+76\bar{\bu}_{j}+67\bar{\bu}_{j+1}-9h\Delta\bu_{j-1}-21h\Delta\bu_{j+1}).
\end{equation}
Thus the linear weights of the three stencils are
\begin{equation}\label{eq:linear-weight}
  \gm^{(-1)} = \dfr {9}{80}, \ \ \gm^{(0)} = \dfr{29}{80}, \ \ \gm^{(1)} = \dfr{21}{40},
\end{equation}
which ensure 
\begin{equation*}
  \tilde\bu_{j+\frac 12,-} = \d\sum_{r=-1}^1\gm^{(r)}\bu^{(r)}_{j+\frac 12, -}.
\end{equation*}
The smoothness indicators are defined by
\begin{equation}\label{eq:SI-def}
  \beta^{(r)} = \d\sum_{l=1}^2 \int_{I_j} h^{2l-1}\left(\dfr{d^l}{dx^l}p^{(r)}(x)\right)^2dx, \ \ r=-1,0,1,
\end{equation}
in the same way as in the WENO reconstructions where $ p^{(r)}(x) $ is the interpolation polynomial on stencil $ S^{(r)} $. Their explicit expressions are
\begin{equation}\label{eq:SI-res}
\begin{array}{l}
  \beta^{(-1)} = (-2\bar{\bu}_{j-1}+2\bar{\bu}_{j}-h\Delta\bu_{j-1})^2 + \dfr{13}{3}(-\bar{\bu}_{j-1}+\bar{\bu}_{j}-h\Delta\bu_{j-1})^2,\\[3mm]
  \beta^{(0)} = \dfr 14 (-\bar{\bu}_{j-1}+\bar{\bu}_{j+1})^2 + \dfr{13}{12}(-\bar{\bu}_{j-1}+2\bar{\bu}_{j}-\bar{\bu}_{j+1})^2,\\[3mm]
  \beta^{(1)} = (2\bar{\bu}_{j+1}-2\bar{\bu}_{j}-h\Delta\bu_{j+1})^2 + \dfr{13}{3}(\bar{\bu}_{j+1}-\bar{\bu}_{j}-h\Delta\bu_{j+1})^2.
\end{array}
\end{equation}
 Then we compute the nonlinear weights in the same way as the WENO-Z method does
\begin{equation}\label{eq:nonlinear-weight}
  \omega^\text{z}_r = \dfr{\al^\text{z}_r}{\sum_{l}\al_l}, \ \ \ \al^\text{z}_r = \gm^{(r)}(1+\dfr{\tau^\text{z}}{\beta^{(r)}+\varepsilon}), \ \ \ r = -1,0,1,
\end{equation}
where $ \varepsilon $ is a small parameter in order to avoid a zero denominator and $ \tau^\text{z} = |\beta^{(1)}-\beta^{(-1)}| $. Finally we have
\begin{equation}\label{eq:final}
  \bu_{j+\frac 12, -} = \d\sum_{r=-1}^1 \omega^\text{z}_r \bu^{(r)}_{j+\frac 12, -}.
\end{equation}
The right interface value $ \bu_{j-\frac 12, +} $ can be reconstructed in a similar way by mirroring the above procedure with respect to $ x_j = \dfr 12(x_{j-\frac 12} + x_{j+\frac 12}) $.\\

Since the GRP solver has to use the spatial derivative $(\pt \bu/\pt x)_{j+\frac 12,\pm}$,  we approximate them using the interpolation, 
\begin{equation}\label{eq:ux-lagrange}
  \Big(\dfr{\pt \bu}{\pt x}\Big)_{j+\frac 12,\pm} := \dfr{1}{12h}\left(\bar\bu_{j-1}-15\bar\bu_{j}+15\bar\bu_{j+1}-\bar\bu_{j+2}\right). 
\end{equation}
The practical simulations later on indicate that the WENO-type stencil selection procedure  in \eqref{eq:ux-lagrange}  can be avoided and the similar observation can be found in \cite{lw-weno, lw-weno-alt}.\\

In \cite{Qiu-11}, the approaches for deriving $\De \bu_j$ were proposed both for the DG method and the finite volume method. In the current study, we use the GRP solver to obtain it, without extra manipulation of governing equations in Section \ref{sec-moment}.  \\

\section{ Construction of gradients  based on the GRP solver} \label{sec-moment}

It is already presented in the original  GRP scheme \cite{Ben-Artzi-01, Li-1} using interface values for the gradient reconstruction. Now we want to apply the idea for the HWENO reconstruction procedure.\\

\subsection{One-dimensional case}\label{sec:moment-1d} First we discuss the one-dimensional case.  Over the computational  cell $I_j$, we regard $\De \bu_j$ as the average of the corresponding spatial derivative, 
\begin{equation}\label{eq:N-L}
  \De \bu_j = \dfr 1h \int_{I_j}\dfr{\pt \bu}{\pt x}(x,t)dx = \dfr 1h \left(\bu(x_{j+\frac 12}, t) - \bu(x_{j-\frac 12}, t)\right),
\end{equation}
where the second equality is the Newton-Leibniz formula. Assume that the GRP values $(\bu_{j+\frac 12}^n, (\pt \bu/\pt t)_{j+\frac 12}^n)$ and $(\bu_{j+\frac 12}^{n+\frac 12}, (\pt \bu/\pt t)_{j+\frac 12}^{n+\frac 12})$ are available around each grid point $x=x_{j+\frac 12}$. Then we obtain the interface values for any time $t\in (t_n,t_{n+1})$. In particular, we have
\begin{equation}
\d \hat{\bu}^{n+\frac 12}_{j+\frac 12} = \bu_{j+\frac 12}^n + \frac {k}{2}\Big(\dfr{\pt\bu}{\pt t}\Big)_{j+\frac 12}^n, \ \ \ \hat{\bu}^{n+1}_{j+\frac 12} = \bu_{j+\frac 12}^n + k\Big(\dfr{\pt\bu}{\pt t}\Big)_{j+\frac 12}^{n+\frac 12}.
\label{eq:cell-boundary}
\end{equation}
It turns out that 
\begin{equation}\label{eq:diff-inter}
  \De \bu^{n+\frac 12}_j = \dfr 1h \left(\hat{\bu}^{n+\frac 12}_{j+\frac 12} - \hat{\bu}^{n+\frac 12}_{j-\frac 12}\right), \ \ \   \De \bu^{n+1}_j = \dfr 1h \left(\hat{\bu}^{n+1}_{j+\frac 12} - \hat{\bu}^{n+1}_{j-\frac 12}\right).
\end{equation}
These values, together with the solution averages $ \bar{\bu}_j^{n+\frac 12} $ and $ \bar{\bu}_j^{n+1} $, are   used in the HWENO reconstruction  at the intermediate stage $ t^{n+\frac 12} $ and the final time stage $ t^{n+1} $, respectively.\\

Now we analyze that such a construction  has the desired accuracy. This is not obvious because the interface values $ \hat{\bu}_{j+\frac 12}^{n+\frac 12} $ and $ \hat{\bu}_{j+\frac 12}^{n+1} $ bear truncation errors of  lower orders compared with those of the cell averages. That is,   $ \hat{\bu}_{j+\frac 12}^{n+\frac 12} $ is first order accurate since it is evaluated by the forward Euler evolution in \eqref{eq:1d-inter} and $ \hat{\bu}_{j+\frac 12}^{n+1} $ is second order accurate with the mid-point rule in \eqref{eq:1d-final}. One might wonder whether the reconstruction can still achieve the desired order of accuracy.

Recall that the  fourth-order accurate numerical  flux is defined as
\begin{equation}\label{eq:flux-4th}
\bga{rl}
  k \bbf_{j+\frac 12}^{4th}& = k \bbf(\bu_{j+\frac 12}^n) +\dfr{k^2}{2} \Big[\dfr 13 \dfr{\pt \bbf}{\pt \bu}(\bu^{n}_{j+\frac 12})\Big(\dfr{\pt \bu}{\pt t}\Big)^{n}_{j+\frac 12} +\dfr 23 \dfr{\pt \bbf}{\pt \bu}(\bu^{n+\frac 12}_{j+\frac 12})\Big(\dfr{\pt \bu}{\pt t}\Big)^{n+\frac 12}_{j+\frac 12}\Big]\\[3mm]
  & = k \bbf(\bu_{j+\frac 12}^n) +\dfr{k^2}{2} \Big[- \dfr 13 \Big(\dfr{\pt \bbf}{\pt \bu}(\bu^{n}_{j+\frac 12})\Big)^2\Big(\dfr{\pt \bu}{\pt x}\Big)^{n}_{j+\frac 12} - \dfr 23 \Big(\dfr{\pt \bbf}{\pt \bu}(\bu^{n+\frac 12}_{j+\frac 12})\Big)^2\Big(\dfr{\pt \bu}{\pt x}\Big)^{n+\frac 12}_{j+\frac 12}\Big] \\[3mm]
  & =\d  \int_{t^n}^{t^{n+1}}\bbf(\bu(x_{j+\frac 12},t))dt +\mathcal{O}(k^5),
\eda
\end{equation}
 where $ (\bu^{n}_{j+\frac 12}, (\pt\bu/\pt x)^{n}_{j+\frac 12}) $ and $ (\bu^{n+\frac 12}_{j+\frac 12}, (\pt\bu/\pt x)^{n+\frac 12}_{j+\frac 12}) $ are the GRP values. This shows that the tolerance for the errors of $ \bu_{j+\frac 12}^n $ is $ \mathcal{O}(k^4) $. And the tolerance for the errors of $ \bu_{j+\frac 12}^{n+\frac 12} $, $ (\pt\bu/\pt x)_{j+\frac 12}^{n} $ and $ (\pt\bu/\pt x)_{j+\frac 12}^{n+\frac 12} $ is $ \mathcal{O}(k^3) $.\\

The Taylor expansion  for $\hat\bu_{j+\frac 12}^{n+\frac 12}$ gives
\begin{equation}\label{eq:boundary-inter}
 \hat{\bu}^{n+\frac 12}_{j+\frac 12} = \bu_{j+\frac 12}^n + \dfr k2 \left(\dfr{\pt\bu}{\pt t}\right)_{j+\frac 12}^n = \bu(x_{j+\frac 12},t^{n+\frac 12}) - \dfr {k^2}8 \left.\dfr{\pt^2\bu}{\pt t^2}\right|_{j+\frac 12}^n + \mathcal{O}(k^3).
\end{equation}
By the definition of  $\De \bu_j^{n+\frac 12}$ in \eqref{eq:diff-inter}, we have
\begin{equation}\label{eq:diff-inter-err}
\begin{array}{l}
  h\De\bu^{n+\frac 12}_j = \hat{\bu}^{n+\frac 12}_{j+\frac 12} - \hat{\bu}^{n+\frac 12}_{j-\frac 12}\\
\ \ \ \ \ \ \ \ \ \ \ = \bu(x_{j+\frac 12},t^{n+\frac 12}) - \bu(x_{j-\frac 12},t^{n+\frac 12}) - \dfr {k^2}{8} \left(\left.\dfr{\pt^2\bu}{\pt t^2}\right|_{j+\frac 12}^n - \left.\dfr{\pt^2\bu}{\pt t^2}\right|_{j-\frac 12}^n\right) + \mathcal{O}(k^4)\\
\ \ \ \ \ \ \ \ \ \ \ = \d\int_{I_j}\dfr{\pt\bu}{\pt x}(x,t^{n+\frac 12})dx - \dfr {k^2h}{8} \left.\dfr{\pt^3\bu}{\pt t^2\pt x}\right|_{j-\frac 12}^n+ \mathcal{O}(k^4).
\end{array}
\end{equation}
 Since $ h $ and $ k $ are proportional thanks  to the CFL condition, $ h\De\bu^{n+\frac 12}_j $ bears the truncation  error  $ \mathcal{O}(k^3) $. Therefore  $ \bu_{j+\frac 12,\pm}^{n+\frac 12} $ is approximated with error  $ \mathcal{O}(k^3) $ in \eqref{eq:s-total}. Finally, since the Riemann solution $ \bu_{j+\frac 12}^{n+\frac 12} $ is calculated from the  data $ \bu_{j+\frac 12,\pm}^{n+\frac 12} $, we conclude that it has the error  $ \mathcal{O}(k^3) $.\\

At the time step $ t^{n+1} $,  it is obvious that the cell average $\bar\bu_j^{n+1}$  has the error $ \mathcal{O}(k^5) $, as shown in \cite{Du-Li-1}. As for the cell boundary values  $\hat\bu_{j+\frac 12}^{n+1}$, they are second order accurate thanks to  the mid-point rule, 
\begin{equation}
  \d \hat{\bu}_{j+\frac 12}^{n+1} = \bu_{j+\frac 12}^n + k\left(\dfr{\pt\bu}{\pt t}\right)_{j+\frac 12}^{n+\frac 12} = \bu(x_{j+\frac 12}, t^{n+1} ) - \dfr{k^3}{24}\left.\dfr{\pt^3\bu}{\pt t^3}\right|_{j+\frac 12}^{n} + \mathcal{O}(k^4), 
\end{equation}
which further gives 
\begin{equation}\label{eq:diff-final-err}
\begin{array}{l}
  h\De\bu^{n+1}_j = \d\int_{I_j}\dfr{\pt\bu}{\pt x}(x,t^{n+1})dx - \dfr {k^3h}{24} \left.\dfr{\pt^4\bu}{\pt t^3\pt x}\right|_{j-\frac 12}^n + \mathcal{O}(k^5).
\end{array}
\end{equation}
With the same arguement,  we can show that  $ h\De\bu^{n+1}_j $ bears an error  $ \mathcal{O}(k^4) $. Therefore the Riemann solution $ \bu_{j+\frac 12}^{n+1} $ bears $ \mathcal{O}(k^4) $.\\

With \eqref{eq:ux-lagrange}, we can show that $ (\pt\bu/\pt x)_{j+\frac 12,\pm}^{n+\frac 12} $ and $ (\pt\bu/\pt x)_{j+\frac 12,\pm}^{n} $ bear errors of orders $ \mathcal{O}(k^3) $ and $ \mathcal{O}(k^4) $, respectively, which meet the above requirement.\\

\subsection{Two-dimensional cases} In this subsection, we supress the dependence of $ \bu $ on the variable $ t $ to simplify the presentation. The gradient construction  for two-dimensional cases can be treated similarly.  Let $\Om_J$ be a computational cell, with the boundary $L_{J\ell}$, $\ell=1,\cdots, K$.  Then we use  the Gauss theorem to have 
\begin{equation}\label{eq:Gauss}
\nb \bu_J : =\dfr{1}{|\Om_J|} \int_{\Om_J} \nb\bu(x,y)dxdy =\dfr{1}{|\Om_J|} \sum_{\ell=1}^K\int_{L_{J_\ell}} \bu   \bn_{J_\ell}dL,
\end{equation}
where $\bn_{J_\ell}$ is the unit outer normal of $\Om_J$ on $L_{J_\ell}$. Therefore, once we know the interface values on $L_{J_\ell}$, we can approximate the gradient $\nb\bu_J$ as 
\begin{equation}\label{eq:gradient-def}
\nb \bu_J \approx \dfr{1} {|\Om_J|} \sum_{\ell=1}^K \sum_m^M \om_{J_{\ell m}}\hat\bu_{J_{\ell m}}   \bn_{J_{\ell m}} |L_{J_\ell}| =:
\left[
  \begin{array}{c}
    \De_x\bu_J\\
    \De_y\bu_J
  \end{array}
\right],
\end{equation} 
where $\hat\bu_{J_{\ell m}}$ is the interface value at the Gauss point $\bx_{J_{\ell m}}$ on the interface $L_{J_\ell}$ and $\om_{J_{\ell m}}$ is the corresponding Gauss weight.\\

Specified to the uniformly rectangular meshes $\Om_{J}=\Om_{ij}=[x_{i-\frac 12},x_{i+\frac 12}]\times [y_{j-\frac 12},y_{j+\frac 12}]$ for which  $ h_x=x_{i+\frac 12}-x_{i-\frac 12} $ and $ h_y=y_{i+\frac 12}-y_{i-\frac 12} $, there are two Gaussian quadrature points on each boundary of $ \Om_{ij} $ to achieve the fifth order accuracy in space. For example we have $ (x_{i+\frac 12},y_{j_1}) $ and $ (x_{i+\frac 12},y_{j_2}) $ on $ x=x_{i+\frac 12}$ for $y\in [y_{j-\frac 12},y_{j+\frac 12}] $, with $y_{j_1}=\frac{1+\sq{3}}{2}y_{j-\frac 12}+\frac{1-\sq{3}}{2}y_{j+\frac 12}$ and $y_{j_2}=\frac{1-\sq{3}}{2}y_{j-\frac 12}+\frac{1+\sq{3}}{2}y_{j+\frac 12}$. The corresponding weights are taken as $ \om_1=\om_2=\frac 12 $. As for the outer normals, we have
\begin{equation}\label{eq:out-norm}
  \bn_{i\pm\frac 12,j_m} = \left[
    \begin{array}{c}
      \pm 1\\
      0
    \end{array}
  \right], \ \ \  \bn_{i_m,j\pm\frac 12} = \left[
    \begin{array}{c}
      0\\
      \pm 1
    \end{array}
  \right], \ \ \ m=1,2.
\end{equation}
It turns out that, by combining \eqref{eq:gradient-def} and \eqref{eq:out-norm}, the gradient $\nabla \bu_{ij}$ is approximated as 
\begin{equation}\label{eq:delta-x-y}
\bga{l}
  \De_x\bu_{ij} = \dfr{1}{h_x}\left[\dfr 12(\hat\bu_{i+\frac 12,j_1}+\hat\bu_{i+\frac 12,j_2})-\dfr 12(\hat\bu_{i-\frac 12,j_1}+\hat\bu_{i-\frac 12,j_2})\right],\\[3mm]
  \De_y\bu_{ij} = \dfr{1}{h_y}\left[\dfr 12(\hat\bu_{i_1,j+\frac 12}+\hat\bu_{i_2,j+\frac 12})-\dfr 12(\hat\bu_{i_1,j-\frac 12}+\hat\bu_{i_2,j-\frac 12})\right],
\eda
\end{equation}
which are the componentwise expressions of \eqref{eq:gradient-def}  over the rectangular grids.\\

Now  we need the limiting values  of the solution $ \bu $ and its derivatives at each Gaussian quadrature point  on the cell boundaries, i.e.,  $ \bu_{i_m,(j+\frac 12,\pm)} $, $ (\frac{\pt \bu}{\pt y})_{i_m,(j+\frac 12,\pm)} $, $ (\frac{\pt \bu}{\pt x})_{i_m,(j+\frac 12,\pm)} $, $ \bu_{(i+\frac 12,\pm),j_m} $, $ (\frac{\pt \bu}{\pt x})_{(i+\frac 12,\pm),j_m} $ and  $ (\frac{\pt \bu}{\pt y})_{(i+\frac 12,\pm),j_m} $ for $ m=1,2 $. We adopt the dimension-by-dimension strategy conventionally  used over rectangular grids \cite{Barth} to interpolate them.

With $ (\bar\bu_{ij},\De_x\bu_{ij}) $, we implement the HWENO reconstruction on $\bu$ in the $x$-direction to obtain the line average of $\bu$ over $ y\in[{y_{j-\frac 12}},{y_{j+\frac 12}}] $ at $ x_{i_1} $ and $ x_{i_2} $, i.e.,
\begin{equation}\label{eq:ux-aver}
  \overline{\bu(x_{i_m},\cdot)}_j := \dfr 1{h_y}\int_{y_{j-\frac 12}}^{y_{j+\frac 12}}\bu(x_{i_m},y)dy, \ \ m=1,2.
\end{equation}
Different from the procedure  for one-dimensional case in Subsection \ref{sec:HWENO}, we interpolate at the Gaussian quadrature points $ x_{i_1} $ and $ x_{i_2} $ instead of the boundaries. Thus  the formulae \eqref{eq:s-0} to \eqref{eq:linear-weight} are modified correspondingly. Furthermore, we have
\begin{equation}\label{eq:ux-diff}
  \De\bu(x_{i_m},\cdot)_j = \dfr{1}{h_y}(\hat\bu_{i_m,j+\frac 12}-\hat\bu_{i_m,j-\frac 12}) \ \ m=1,2.
\end{equation}
Finally, we implement the HWENO reconstruction for  $ \bu(x_{i_m},\cdot) $ with $ (\overline{\bu(x_{i_m},\cdot)}_j,(\De_y\bu(x_{i_m},\cdot))_j) $ defined in \eqref{eq:ux-aver} and \eqref{eq:ux-diff} to obtain $ \bu(x_{i_m},y_{j+\frac 12,\pm}) $ and $ \frac{\pt \bu}{\pt y}(x_{i_m},y_{j+\frac 12,\pm}) $,  which are just $ \bu_{i_m,(j+\frac 12,\pm)} $ and $ (\frac{\pt \bu}{\pt y})_{i_m,(j+\frac 12,\pm)} $.

As for the tangential derivatives $ (\frac{\pt \bu}{\pt x})_{i_m,(j+\frac 12,\pm)} $, implement the HWENO reconstruction on $\bu$ with $ (\bar\bu_{ij},\De_y\bu_{ij}) $ in the $y$-direction to obtain the line average of $\bu$ over $ x\in[{x_{i-\frac 12}},{x_{i+\frac 12}}] $ at $y_{j+\frac 12,\pm}$, i.e.,
\begin{equation*}
  \overline{\bu(\cdot,y_{j+\frac 12,\pm})}_i := \d \frac 1{h_x} \int_{i-\frac 12}^{i+\frac 12}\bu(x,y_{j+\frac 12,\pm})dx.
\end{equation*}
With $\overline{\bu(\cdot,y_{j+\frac 12,\pm})}_i$ already obtained above, we interpolate the derivatives of $ \bu(\cdot,y_{j+\frac 12,\pm}) $ at $x_{i_1}$ as
\begin{equation}
\bga{rl}
 & \dfr{\pt\bu}{\pt x}(x_{i_1},y_{j+\frac 12,\pm})\\[3mm]
=&\dfr{1}{108h_x}\left[(9+2\sq{3})\overline{\bu(\cdot,y_{j+\frac 12,\pm})}_{i-2}-(72+26\sq{3})\overline{\bu(\cdot,y_{j+\frac 12,\pm})}_{i-1}\right.\\[3mm]
&\left.+48\sq{3}~\overline{\bu(\cdot,y_{j+\frac 12,\pm})}_{i} + (72-26\sq{3})\overline{\bu(\cdot,y_{j+\frac 12,\pm})}_{i+1} - (9-2\sq{3})\overline{\bu(\cdot,y_{j+\frac 12,\pm})}_{i+2}\right],
\eda
\label{eq:deri-gauss}
\end{equation}
which is just $ (\frac{\pt \bu}{\pt x})_{i_1,(j+\frac 12,\pm)} $,  the tangential derivative  of $ \bu $ at the Gaussian quadrature point $ (x_{i_1},y_{j+\frac 12,\pm}) $. We can obtain $ (\frac{\pt \bu}{\pt x})_{i_2,(j+\frac 12,\pm)} $ by simply mirroring \eqref{eq:deri-gauss} with respect to $ (x_i,y_{j+\frac 12,\pm}) $.\\

Similar interpolations are made for $ \bu_{(i+\frac 12,\pm),j_m} $, $ (\frac{\pt \bu}{\pt x})_{(i+\frac 12,\pm),j_m} $ and  $ (\frac{\pt \bu}{\pt y})_{(i+\frac 12,\pm),j_m} $ for $ m=1,2 $. The details are omitted. \\

%\begin{rem}
For  system cases, the characteristic decomposition  in \cite{Barth} is adopted in the current study. All the reconstruction procedures described above are applied for the characteristic variables. The details are omitted here.\\
%\end{rem}

\vspace{2mm}

\section{Numerical Examples}\label{sec:numer}
In this section, we provide  several examples for one- and two-dimensional compressible Euler equations to verify the expected performance of this approach. The Euler equations can be found in any CFD books and we do not write out them here. Each example will be computed with the HWENO interpolation and the WENO interpolation in the  same framework of  two-stage fourth-order time discretization based the GRP solver \cite{Du-Li-1}. Both the HWENO and the WENO reconstruction use the nonlinear weights in formulae \eqref{eq:nonlinear-weight}. The resulting schemes are denoted as GRP4-HWENO5 and GRP4-WENO5, respectively. We emphasize again that the only difference is the data reconstruction. The CFL number $ 0.6 $ is used for all the computations except in the first example.\\ 
 
\n{\bf Example 1. Smooth initial value problem.} 
We check the numerical results for a one-dimensional smooth initial value problem of the Euler equations with the initial data
\begin{equation}\label{eq:smooth}
     \rho(x,0)   =     1+0.2\text{sin}(\pi x), \ \ \ v(x,0)=1, \ \  p(x,0)=1,
\end{equation}
to verify the numerical accuracy of the present approach where $ \rho $ is the density, $ v $ is the velocity, $ p $ is the pressure and $ \bu=[\rho,v,p]^\top $. The periodic boundary conditions are used. The exact solution at time $ t $ is just a shift of the initial condition, i.e.
\begin{equation}\label{iv:smooth}
\bu(x,t)=\bu(x-t,0).
\end{equation}
Here we set the CFL number to be $0.1$ to show the numerical order of the spatial reconstructions. The errors shown in Table \ref{tab:Euler_smooth} are those of the cell averages of the density $ \rho $ at time $ t=10 $.

Both reconstruction approaches achieve the designed numerical order while the errors of the scheme GRP4-HWENO5 is smaller than those of GRP4-WNEO5. And we can see that the the CPU time cost by both schemes are almost the same which verifies our claim that no additional efforts are made to obtain the first moment of the solution.% Their convergence rates are shown in Figure \ref{fig:error-Euler}  

\begin{table*}[!htbp]
  \centering
  \caption[small]{The $L_1$, $L_\iy$ errors of the density and numerical orders for the smooth initial value problem in Example 1. The results are shown at time $ t=10 $.}
\small
  \begin{tabular}{|r|r|l|l|l|l|r|l|l|l|l|}
    \hline
$m$&     \multicolumn{5}{|c|}{GRP4-WENO5}      &       \multicolumn{5}{|c|}{GRP4-HWENO5}    \\\cline{2-11}
   &CPU time (s)&$ L_1 $ error&Order&$ L_\infty $ error&Order&CPU time (s)&$ L_1 $ error &Order&$ L_\infty $ error &Order\\\hline
40  &    0.31  &  6.25e-6  &  4.98  &  9.83e-6  &  4.99  &  0.41  &  1.50e-6  &  4.99  &  2.36e-6  &  5.02\\
80  &    1.25  &  1.96e-7  &  4.99  &  3.08e-7  &  5.00  &  1.46  &  4.69e-8  &  5.00  &  7.37e-8  &  5.00\\
160  &   5.14  &  6.13e-9  &  5.00  &  9.64e-9  &  5.00  &  4.95  &  1.47e-9  &  5.00  &  2.30e-9  &  5.00\\
320  &   19.77 &  1.92e-10  &  5.00  &  3.01e-10  &  5.00  &  19.61 &  4.59e-11  &  5.00  &  7.21e-11  &  5.00\\
640  &   122.38 & 5.99e-12  &  5.00  &  9.47e-12  &  4.99  &  117.23 & 1.43e-12  &  5.00  &  2.36e-12  &  4.93\\\hline
  \end{tabular}
  \label{tab:Euler_smooth}
\end{table*}
\vspace{0.2cm}
\iffalse
\begin{figure}[!htb]
\centering
 \includegraphics[width=.6\textwidth]{./error-sin-1d.eps}
 \caption[small]{The $L_1$ and $L_\iy$ convergence rates of GRP4-WENO5 and GRP4-HWENO5 for the Euler equations in Example 1.}
\label{fig:error-Euler}
\end{figure}
\fi
\vspace{0.2cm} 

\begin{figure}[!htb]
 \centering
 \includegraphics[width=.9\textwidth]{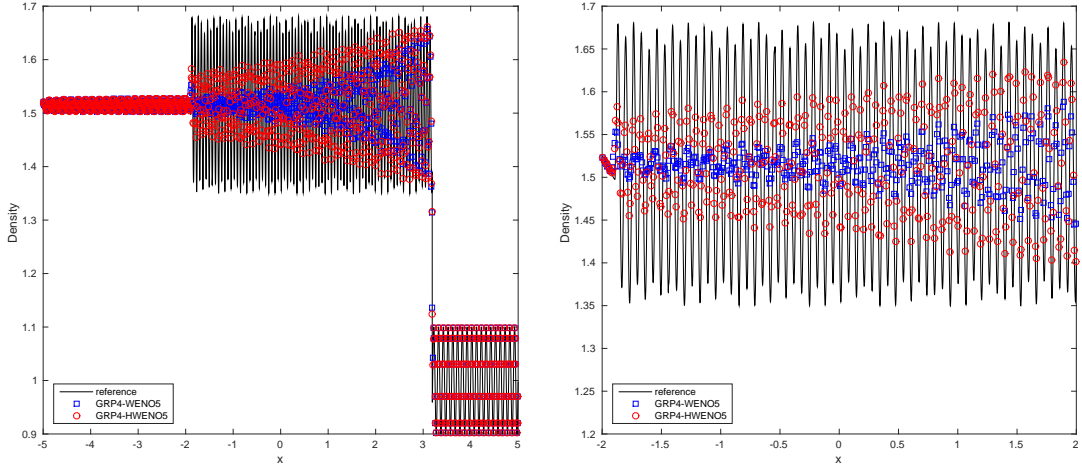}
  \caption[small]{The comparison of the density profile (left) and its local enlargement (right) for the Titarev-Toro problem in Example 2. The schemes used are GRP4-WENO5 (squares) and GRP4-HWENO5 (circles) with $1000$ cells. The solid lines are the reference solution.}
 \label{fig:Toro}
\end{figure}

 \n{\bf Example 2.  The Titarev-Toro problem.}   This example was proposed in \cite{Toro-14} as an extension of the Shu-Osher problem \cite{shu-osher}. The initial data is taken as
 \begin{equation}
(\rho,v,p)(x,0)=\left\{ \begin{array}{ll}
(1.515695, 0.523346, 1.805), \ \ \ \ &\mbox{for } x<-4.5,  \\[3mm] 
  (1 + 0.1 \sin(20 \pi x), 0, 1),&\mbox{for } x\geq -4.5. 
 \end{array}
 \right.
  \end{equation} 
The output time is $ t=5 $ and the numerical solutions computed with $ 1000 $ cells are shown in Figure \ref{fig:Toro}. The reference solution is computed with $10000$ cells. The scheme GRP4-HWENO5 can catch the peaks and the troughs in the solution better.
  
\vspace{0.2cm}

 \n{\bf Example 3. Large pressure ratio problem. }  The large pressure ratio problem is a Riemann problem first presented in \cite{Tang-Liu}. In this problem, initially the pressure and density ratio between the two neighboring states are very high. The initial data is $(\rho,v,p)=(10000,0,10000)$ for $0\leq x<0.3$ and $(\rho,v,p)=(1,0,1)$ for $0.3\leq x \leq 1.0$. The boundary condition is dealt with in the same with that in the standard Riemann problem. The output time is $ t = 0.12 $.\\

\begin{figure}[htp]
\centering
\includegraphics[width=.9\textwidth]{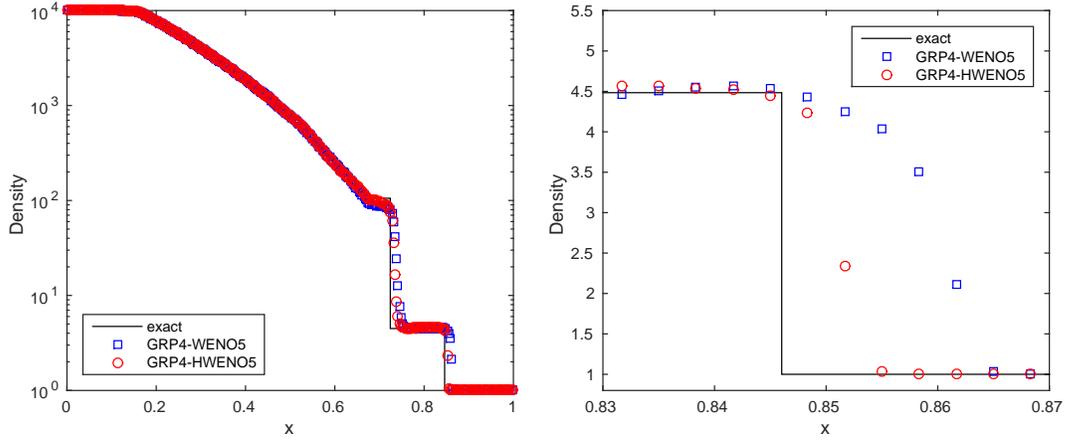}
\caption{The comparison of the density profile (left) and its local enlargement (right) for the large pressure ration problem in Example 3. The schemes GRP4-WENO5 (squares) and GRP4-HWENO5 (circles) are performed with $ 300 $ cells. The solid lines refer to  the exact solution.}
\label{fig:ratio}
\end{figure}

An extremely strong rarefaction wave  forms and it  significantly affects the shock location in the numerical solution. Two factors determine the ability of a numerical scheme to properly capture the position of the shock. The first one is whether the thermodynamical effect is included in the numerical fluxes properly \cite{thermo-GRP}. The second one is the numerical dissipation of the schemes. Figure \ref{fig:ratio} shows the density profile of the numerical solutions  simulated by GRP4-WENO5 and GRP4-HWENO5 with $ 300 $ cells. We can see that GRP4-HWENO5 behaves better due to   smaller stencils and less numerical dissipations of the resulting scheme.

\vspace{0.2cm}

\begin{figure}[htp]
\centering
\includegraphics[width=.8\textwidth]{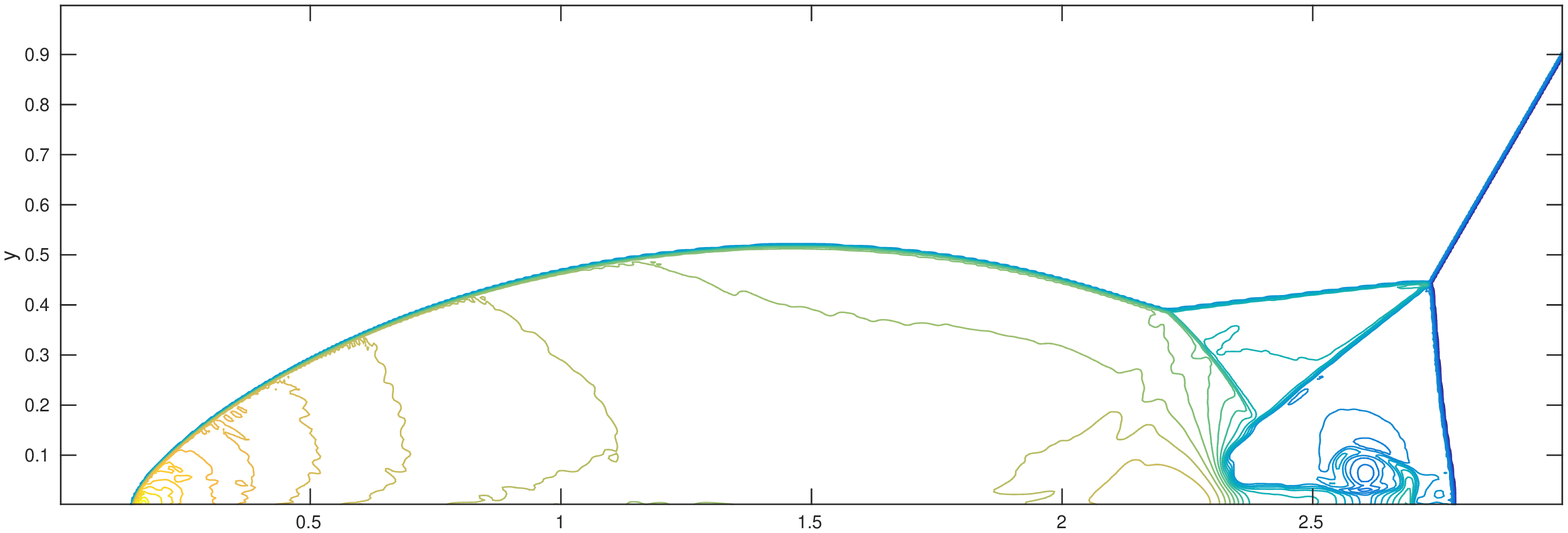}
\includegraphics[width=.8\textwidth]{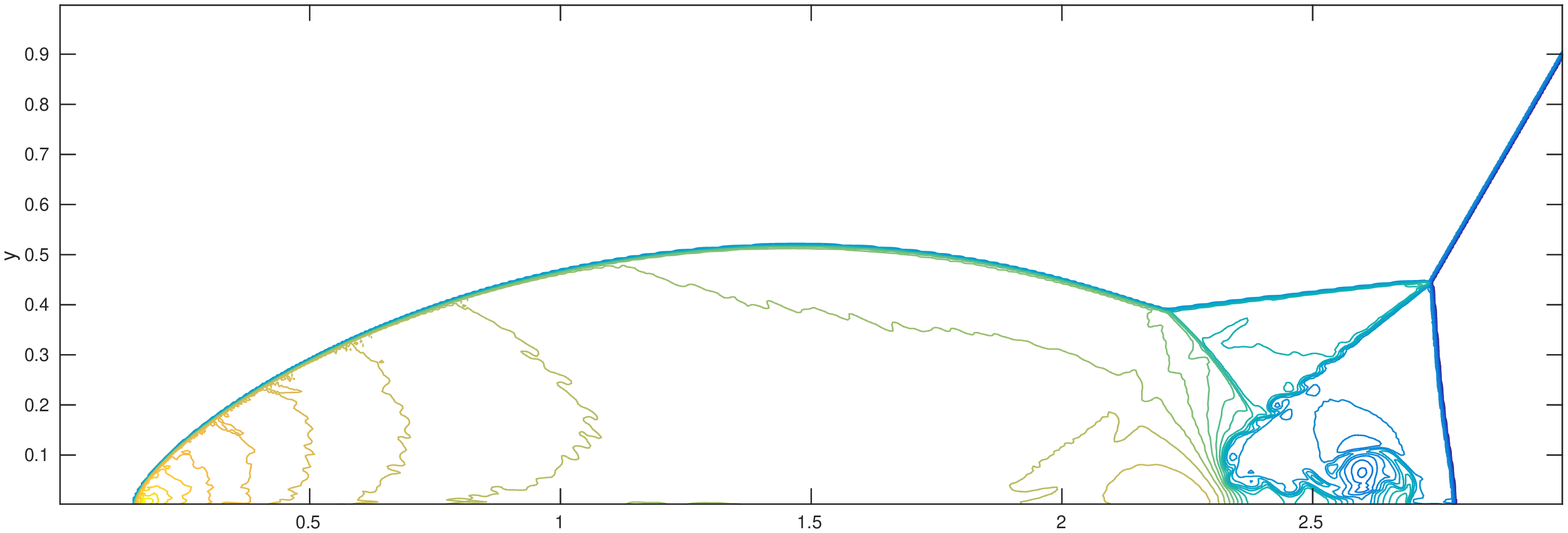}
\caption{The density contours of the double Mach reflection problem in Example 4  by GRP4-WENO5 (upper) and the GRP4-HWENO5 (lower) with $ 960\times240 $ cells.}
\label{fig:double-mach}
\end{figure}

\n{\bf Example 4. The double Mach reflection problem.} This is a standard test problem to display the performance of high resolution schemes. The computational domain for this problem is $ [0,4] \times [0,1] $, and $ [0,3] \times [0,1]$ is shown. The reflective wall lies at the bottom of the computaional domain starting from $ x=\frac{1}{6} $. Initially a right-moving Mach 10 shock is positioned at $ x=\frac{1}{6} ,y=0$ and makes a $ \frac{\pi}{3} $ angle with the $x$-axis. The results are shown in Figure \ref{fig:double-mach} from which we can see that the numerical result  by  GRP4-HWENO5 can resolve more structures along the slip line than that by GRP4-WENO5.

\vspace{0.2cm}

\n{\bf Example 5.  Two-dimensional Riemann problems.} We provide an example of two-dimensional Riemann problem taken from \cite{Han} involving the interactions of vortex sheets with rarefaction waves. The computation is implemented over the domain $[0,1]\times [0,1]$.

\begin{equation}
(\rho,u,v,p)(x,y,0) =\left\{
\begin{array}{ll}
 (1, 0.1, 0.1, 1), & 0.5<x<1,0.5<y<1,\\
 (0.5197, -0.6259, 0.1, 0.4), & 0<x<0.5, 0.5<y<1,\\
 (0.8, 0.1, 0.1, 0.4),&0<x<0.5,0<y<0.5,\\
 (0.5197, 0.1, -0.6259, 0.4), &0.5<x<1,0<y<0.5.
\end{array}
\right.
\end{equation}
The output time is $0.3$.

\begin{figure}[htp]
\centering
\includegraphics[width=.8\textwidth]{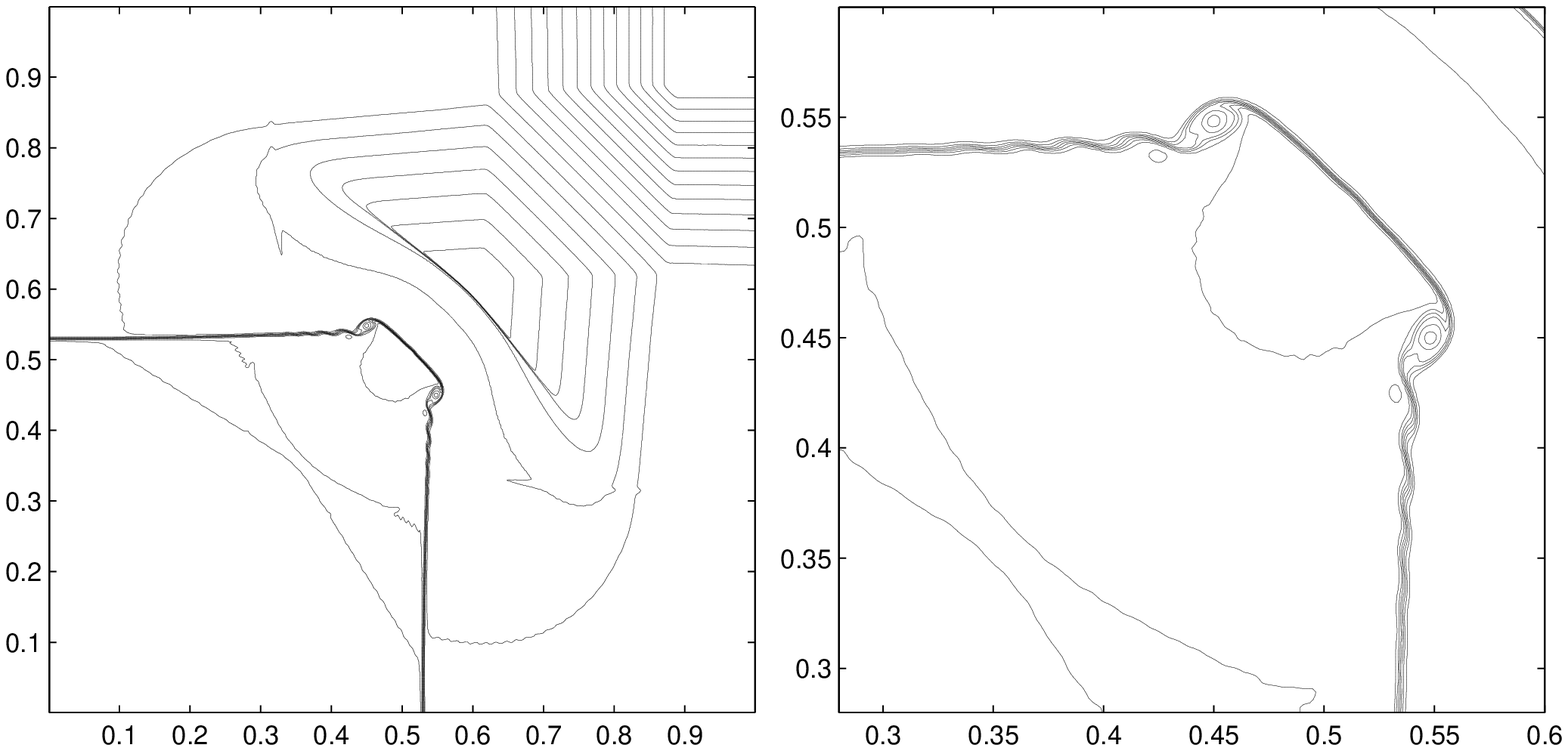}
\includegraphics[width=.8\textwidth]{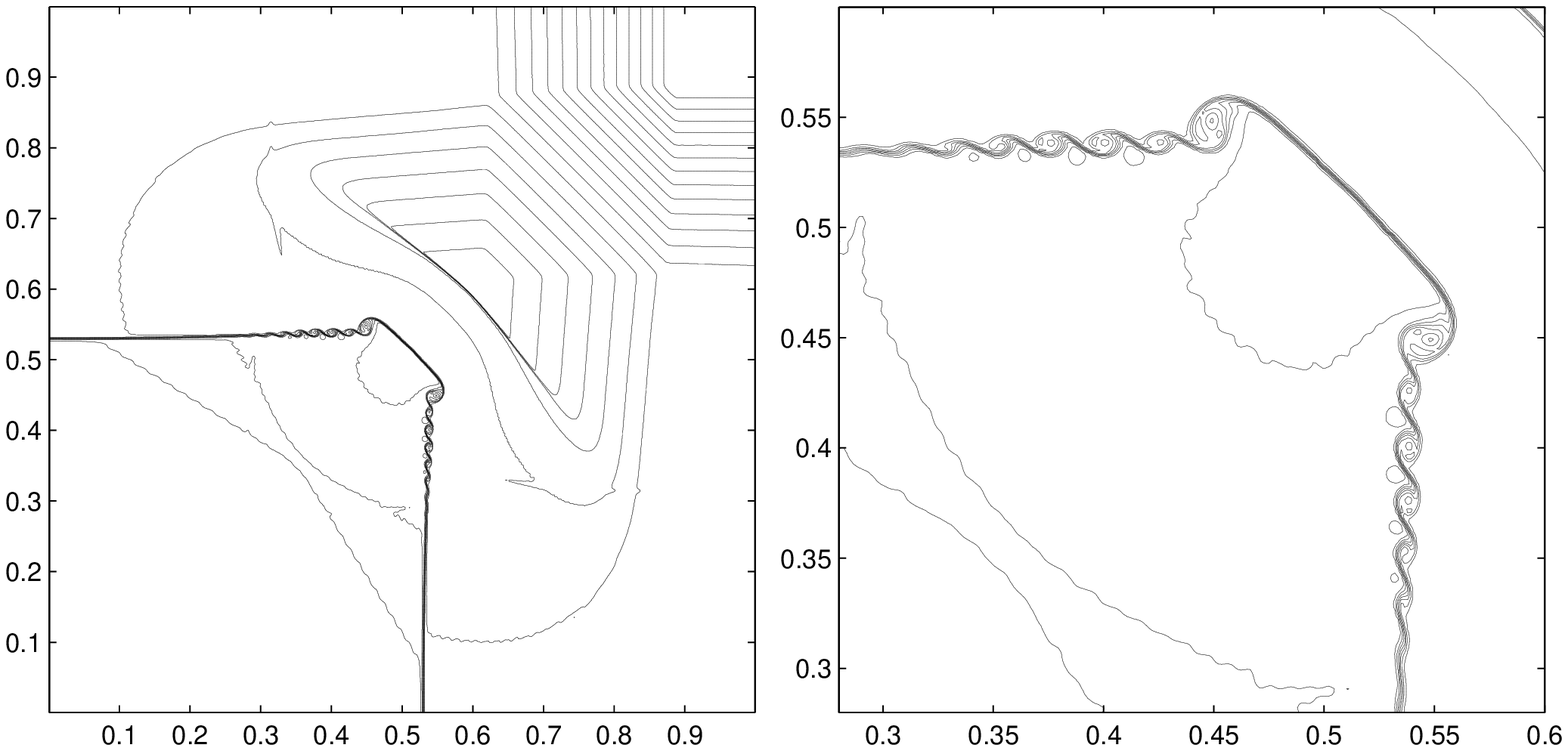}
\caption{The density contours of the 2-D Riemann problem in Example 5 computed with the schemes GRP4-WENO5 (upper) and GRP4-HWENO5 (lower), respectively. $ 700\times700 $ cells are used.}
\label{fig:Riemann-a}
\end{figure}

The contours of the density and their local enlargements are shown in Figures \ref{fig:Riemann-a}. We can see that the scheme with the HWENO reconstruction can resolve more small structures along the vortex sheet.

\vspace{2mm}

\section{Discussion}

In this paper we developed a new HWENO reconstruction method just over structural (rectangular) meshes. This can be extended over unstructured meshes, but we think that the technicality may be nontrivial, which remains for future study. \\

A subtlety also lies in the reconstruction of spatial derivatives close to  interfaces  at the intermediate stage.  If we would use the formulae below to  interpolate the first-order derivatives as,
\begin{equation}\label{eq:ux-hermite}
\bga{l}
  \Big(\dfr{\pt \bu}{\pt x}\Big)_{j-\frac 12,+} := \dfr{1}{8h}\left(-13\bar\bu_{j-1}+16\bar\bu_{j}-3\bar\bu_{j+1} -3h\De\bu_{j-1}+h\De\bu_{j+1}\right),\\[2mm]
  \Big(\dfr{\pt \bu}{\pt x}\Big)_{j+\frac 12,-} := \dfr{1}{8h}\left(3\bar\bu_{j-1}-16\bar\bu_{j}+13\bar\bu_{j+1} + h\De\bu_{j-1}-3h\De\bu_{j+1}\right),
\eda
\end{equation}
the same analysis performed in Subsection \ref{sec:moment-1d} shows  that  the  errors of $ \Big(\dfr{\pt \bu}{\pt x}\Big)_{j+\frac 12,\pm} $ would be  $ \mathcal{O}(k^2) $ so that 
\begin{equation*}\label{eq:defect}
  \bbf_{j+\frac 12}^{4th} =
   \d\dfr{1}{k} \int_{t^n}^{t^{n+1}}\bbf(\bu(x_{j+\frac 12},t))dt+\mathcal{O}(k^3). 
\end{equation*}
This makes the resulting scheme third-order accurate. \\

There is a remedy here. We could use the third order GRP solver in \cite{Qian-Li} , which is relatively complicated and we want to use an alternative method. Instead of the formulae in \eqref{eq:ux-hermite}  to interpolate $ \Big(\dfr{\pt\bu}{\pt x}\Big)^{n+\frac 12}_{j+\frac 12,\pm} $, we adopt 
\begin{equation}\tag{\ref{eq:ux-lagrange}}
  \Big(\dfr{\pt \bu}{\pt x}\Big)_{j+\frac 12,\pm} := \dfr{1}{12h}\left(\bar\bu_{j-1}-15\bar\bu_{j}+15\bar\bu_{j+1}-\bar\bu_{j+2}\right)
\end{equation} 
in practice. Although the stencil becomes larger, the numerical results in Section \ref{sec:numer} show that such a choice has a good practical effect. \\

\vspace{2mm}

\vspace{2mm}

\bibliographystyle{plain}

\vspace{0.5cm}
\end{document}